\documentclass[11pt, oneside]{article}
\usepackage{amsmath,amsfonts,amssymb,amscd,theorem}
\date{}

\newtheorem{thm}{Theorem}
\newtheorem{prop}{Proposition}
\newtheorem{lem}{Lemma}
\newtheorem{cor}{Corollary}
\newtheorem{rem}{Remark}

\begin{document}

\title{Cohomologies of Harmonic Bundles \\ on Quasi-Compact K\"ahler manifolds}

\author{J\"urgen Jost \footnote{Max-Planck Institute for
Mathematics in the Sciences, Leipzig, Germany}, Yi-Hu Yang
\footnote{Department of Mathematics, Tongji University, Shanghai,
China.}
\thanks{The author supported partially by NSF of China (No.
10471105, 10771160)}, and Kang Zuo \footnote{Department of
Mathematics, Mainz University, Mainz, Germany}}

\maketitle

\section{Introduction and the statements of the problems}

Let $\mathbb{V}$ be a complex flat vector bundle of rank $n$ on a
complex manifold $M^m$; equivalently, this corresponds to a linear
representation
$$
\rho: \pi_1(M)\to GL(n, \mathbb{C}).
$$
\\
If $M$ is compact K\"ahlerian and $\rho$ is {\it semi-simple},
that is that the Zariski closure of the image of $\rho$ in $GL(n,
\mathbb{C})$ is semi-simple {\footnote{Equivalently, call $\rho$
semi-simple if for any boundary component $\Sigma$ of
$\mathcal{P}_n:=GL(n, \mathbb{C})/U(n)$, there exists an element
$\gamma\in\pi_1(X)$ satisfying
$\Sigma\cap\rho(\gamma)(\Sigma)=\emptyset$.}}, by means of a
result of Donaldson (in the case of Riemann surfaces, \cite{do})
and Corlette (in the higher dimensional case, \cite{co}), there
exists a unique harmonic metric $u$ on $\mathbb{V}$, equivalently,
a $\rho$-equivariant harmonic map from the universal covering
$\tilde M$ of $M$ into $\mathcal{P}_n:=GL(n, \mathbb{C})/U(n)$--
the set of positive definite Hermitian symmetric matrices of order
$n$.
\\
\\
Let $D$ be the flat connection of $\mathbb{V}$. Decompose
$D=d'+d''$ into $(1, 0)$-part and $(0, 1)$-part. Suppose
$\delta''$ and $\delta'$ are $(0, 1)$-type operator and $(1,
0)$-type operator respectively, s.t. $d'+\delta''$ and
$\delta'+d''$ preserve the Hermitian metric $u$.
\\
\\
According to Simpson \cite{si, si1}, set
\begin{eqnarray*}
&&\partial=(d'+\delta')/2, ~~~{\overline\partial}=(d''+\delta'')/2,\\
&&\theta=(d'-\delta')/2, ~~~{\overline\theta}=(d''-\delta'')/2.
\end{eqnarray*}
It is easy to check that i) $\partial+\overline\partial$ preserves
the metric $u$; ii) $<\theta x, y>_u=<x, \overline\theta y>_u$;
iii) up to some constant, $\theta=(\partial u)u^{-1}$ and $\theta$
can be considered as a $1$-form of type $(1, 0)$ valued in
${\text{Hom}}(\mathbb{V})$, which we call the {\it Higgs field}.
\\
\\
Put
\begin{eqnarray*}
&&D_u''=\overline\partial +\theta,
~~~D'_u=\partial+\overline\theta,\\
&& D^c_u=D''_u-D'_u.
\end{eqnarray*}
\\
When $M$ is one-dimensional, Hitchin \cite{hi} observed that the
harmonicity of $u$ gives rise to the structure of a Higgs bundle
on $\mathbb{V}$; namely, since $(D_u'')^2=0$, one can consider
$\mathbb{V}$ as a holomorphic vector bundle under the operator
$\overline\partial$, denoted by $E$ later on, and hence $\theta$
is a holomorphic $1$-form valued in ${\text{Hom}}(E)$. In the
higher dimensional case, using Siu's Bochner formula (\cite{siu},
also see \cite{samp}), one can show that $u$, as a map, is
pluri-harmonic \cite{co}, equivalently, this is $(D_u'')^2=0$. So,
in the higher dimensional case, one still has the above structure
of a Higgs bundle with $\theta\wedge\theta=0$ (this is automatical
in the case of Riemann surface). An equivalent but more geometric
observation was actually discovered earlier by Jost-Yau \cite{jy},
where from the harmonic map they derived a fibration which is
holomorphic, namely $\theta$ is holomorphic. Later on, we call
$(E, \theta)$ with $\theta\wedge\theta=0$ the {\it Higgs bundle}
and $D_u''$ the {\it Higgs operator}, which is integrable.
\\
\\
{\bf Definition.} {\it A harmonic bundle on a complex manifold $M$
is a complex flat vector bundle $\mathbb{V}$ over $M$ with a
harmonic metric $u$ satisfying the induced operator $D_u''$ is
integrable, so that one has an induced Higgs bundle $(E, \theta)$
as above.}
\\
\\
For a harmonic bundle on a complex manifold, one can define
various cohomologies.
\\
\\
Let us first give some formal definitions. Let $\mathbb{V}$ be a
harmonic bundle with a harmonic metric $u$. By $\mathbb{V}^D$
denote the locally constant sheaf of flat sections of
$\mathbb{V}$. This sheaf is resolved by de Rham complex of sheves
of $C^\infty$ differential forms with coefficients in
$\mathbb{V}$:
\begin{equation} {\begin{CD}
\mathbb{V}^D\to(\{\mathcal{A}^\centerdot(\mathbb{V}),
D\}=\mathcal{A}^0(\mathbb{V})@>{D}>>
\mathcal{A}^1(\mathbb{V})@>{D}>>\mathcal{A}^2(\mathbb{V})@>{D}>>\cdots)
\end{CD}}
\end{equation}
is a quasi isomorphism of complexes of sheaves. The sheaves of
$C^\infty$ forms are fine, so the sheaf cohomology $H^*(M,
\mathbb{V}^D)$ is naturally isomorphic to the cohomology of the
complex of global sections
\begin{equation} {\begin{CD}
\{A^\centerdot(\mathbb{V}), D\}=A^0(\mathbb{V})@>{D}>>
A^1(\mathbb{V})@>{D}>>A^2(\mathbb{V})@>{D}>>\cdots
\end{CD}}
\end{equation}
Call the cohomology of the complex of global sections the {\it de
Rham cohomology}, denoted by $H^*_{\text{{DR}}}(M, \mathbb{V})$.
\\
\\
For the corresponding Higgs bundle $(E, \theta)$, we also can
define cohomologies as follows. First, we have a complex of free
sheaves of locally holomorphic forms valued in $E$, called the
{\it holomorphic Dolbeault complex}:
\begin{equation} {\begin{CD}
\{\Omega^\centerdot,
\theta\wedge\}=E@>{\theta\wedge}>>\Omega^1(E)@>{\theta\wedge}>>\Omega^2(E)@>{\theta\wedge}>>\cdots
\end{CD}}
\end{equation}
The condition that $\theta\wedge\theta=0$ insures that this is a
complex. Then, define the {\it Dolbeault cohomology} with
cofficients in $E$ to be the hypercohomology of the above complex
\begin{equation*} {\begin{CD}
{\bf H}^*(E@>{\theta\wedge}>>\Omega^1(E)@>{\theta\wedge}>>\cdots),
\end{CD}}
\end{equation*}
denoted by $H^*_{\text{Dol}}(M, E)$.
\\
\\
Furthermore, the complex of sheaves of local $C^\infty$ sections
valued in $E$ with the differential $D''_u$
\begin{equation} {\begin{CD}
\{\mathcal{A}^\centerdot(E), D''_u\}=\mathcal{A}^0(E)@>{D''_u}>>
\mathcal{A}^1(E)@>{D''_u}>>\mathcal{A}^2(E)@>{D''_u}>>\cdots
\end{CD}}
\end{equation}
gives a fine resolution of the holomorphic Dolbeault complex, so
$H^*_{\text{Dol}}(M, E)$ is naturally isomorphic to the cohomology
of the complex of global sections $\{A^\centerdot(E), D''_u\}$.
But, for later convenience, we call the cohomology of the complex
of global sections the {\it Higgs cohomology}, denoted by
$H^*_{\text{Higgs}}(M,E)$; and the complex
$\{\mathcal{A}^\centerdot(E), D''_u\}$ the {\it Higgs complex}.
\\
\\
A (formal) discussion tells us that if $M$ is compact K\"ahler,
all the above cohomologies: the sheaf cohomology $H^*(M,
\mathbb{V}^D)$, the de Rham cohomology $H^*_{\text{{DR}}}(M,
\mathbb{V})$, the Dolbeault cohomology $H^*_{\text{Dol}}(M, E)$,
and the Higgs cohomology $H^*_{\text{Higgs}}(M,E)$, are
isomorphic. The above arguments tell us that one only needs to
show that the de Rham cohomology $H^*_{\text{{DR}}}(M,
\mathbb{V})$ is isomorphic to the Higgs cohomology
$H^*_{\text{Higgs}}(M,E)$. This is just a consequence of the
(local) K\"ahler identity for harmonic bundles and the theory of
harmonic forms; for details, see the end of this section.
\\
\\
When $M$ is non-compact, e.g. $M={\overline M}-D$, $D$ being a
divisor of $\overline M$, the situation gets very complicated; one
does not necessarily have the isomorphism between the de Rham
cohomology $H^*_{\text{{DR}}}(M, \mathbb{V})$ and the Dolbeault
cohomology $H^*_{\text{Dol}}(M, E)$. Roughly speaking, these
cohomologies are too big for study and applications, so that we
have to modify all the above four classes of cohomologies. A good
method for modifying them are to consider certain sub-complexes by
using certain growth condition and the corresponding cohomologies.
\\
\\
In order to modify these cohomologies, let us consider the
following setup, which geometrically is the most interesting. Let
$\overline M$ be a compact K\"ahler manifold with a K\"ahler
metric $\omega_0$, $D=\sum_1^s D_k$ a normal crossing divisor,
denote ${\overline M}-D$ by $M$ and the inclusion map
$M\hookrightarrow \overline M$ by $j$. Call $M$ a {\it
quasi-compact K\"ahler manifold}. One can equip $M$ with a
complete K\"ahler metric $\omega$ which is Poincar\'e-like near
the divisor, namely, if $p\in D$ and ${\bf z}=(z_1, z_2, \cdots,
z_m)$ is a local coordinate at $p$ satisfying that ${\bf z}(p)=(0,
0, \cdots, 0)$ and $\prod_1^{s'}z_k=0$, $s'\le s$, is the local
minimal equation of $D$, then, near $p$
\[
\omega=\sum_{k=1}^{s'}{\frac{\sqrt{-1}dz_k\wedge d{\overline
z}_k}{|z_k|^2|\log|z_k||^2}} + \sqrt{-1}\sum_{s'+1}^m dz_k\wedge
d{\overline z}_k.
\]
Now, suppose that $\mathbb{V}$ is a harmonic bundle over $M$ with
a harmonic metric $u$, $E, \theta$ the corresponding the Higgs
bundle, the Higgs field. Using the metric $\omega$ on $M$ and the
metric $u$, one can define (local) square-integrable forms valued
in $\mathbb{V}$ or $E$. Using (local) square-integrable forms, one
then can define the corresponding sub-complexes on $\overline M$
of (1) and (2) respectively as follows:
\begin{equation} {\begin{CD}
\{\mathcal{A}_{(2)}^\centerdot(\mathbb{V}),
D\}=\mathcal{A}_{(2)}^0(\mathbb{V})@>{D}>>
\mathcal{A}_{(2)}^1(\mathbb{V})@>{D}>>\mathcal{A}_{(2)}^2(\mathbb{V})@>{D}>>\cdots,
\end{CD}}
\end{equation}
\begin{equation} {\begin{CD}
\{A_{(2)}^\centerdot(\mathbb{V}), D\}=A_{(2)}^0(\mathbb{V})@>{D}>>
A_{(2)}^1(\mathbb{V})@>{D}>>A_{(2)}^2(\mathbb{V})@>{D}>>\cdots,
\end{CD}}
\end{equation}
where the sheaves $\mathcal{A}_{(2)}^i(\mathbb{V})$ are defined as
follows:
\\
{\it for an open subset $U$ of $\overline M$,
$\mathcal{A}^i_{(2)}(\mathbb{V})(U)$ is defined as the set of
$\mathbb{V}$-valued $i$-forms $\eta$ on $U\cap M$ with measurable
coefficients and measurable exterior derivative $D\eta$, such that
$\eta$ and $D\eta$ have finite $L^2$-norm};
\\
and $\{A_{(2)}^\centerdot(\mathbb{V}), D\}$ is the complex of
global sections of $\{\mathcal{A}_{(2)}^\centerdot(\mathbb{V}),
D\}$.
\\
\\
Since the sheaves $\mathcal{A}_{(2)}^\centerdot(\mathbb{V})$ are
fine, the hypercohomology of the complex
$\{\mathcal{A}_{(2)}^\centerdot(\mathbb{V}), D\}$ is isomorphic to
the cohomology of the complex $\{A_{(2)}^\centerdot(\mathbb{V}),
D\}$. We call this cohomology the {\it $L^2$-de Rham cohomology},
denoted by $H_{{\text{DR}},(2)}^*(M, \mathbb{V})$;
correspondingly, the complex is called the {\it $L^2$-de Rham
complex.}
\\
\\
On the other hand, as mentioned before, the sheaf cohomology
$H^*(M, \mathbb{V}^D)$ is too big in such a non-compact case,
according to a suggestion of Deligne, one should consider the
(middle) {\it intersection cohomology} \cite{cgm, gm1, gm2} with
coefficients in the direct sheaf $j_*\mathbb{V}$:
$H^*_{\text{int}}(\overline M, j_*\mathbb{V})$. Then, one has the
following
\\
\\
{\bf Conjecture 1.} {\it The intersection cohomology
$H^*_{\text{int}}(\overline M, j_*\mathbb{V})$ is isomorphic to
the $L^2$-de Rham cohomology $H_{{\text{DR}},(2)}^*(M,
\mathbb{V})$.}
\\
\\
Now, we turn to the Dolbeault cohomology. In this case, we need to
pay much more attention to the properties at the divisor of the
representation $\rho$ of $\pi_1(M)$ and the asymptotic behavior of
the harmonic metric $u$; but here, we will not involve in much
more details, instead, we will give some direct statements; in the
next two sections, these will be stated more clearly.
\\
\\
Similar to the $L^2$-de Rham complex, we define the {\it
$L^2$-holomorphic Dolbeault complex} and  the {\it $L^2$-Higgs
complex} as follows:
\begin{equation} {\begin{CD}
\{\Omega_{(2)}^\centerdot,
\theta\wedge\}=E_{(2)}@>{\theta\wedge}>>\Omega_{(2)}^1(E)@>{\theta\wedge}>>
\Omega_{(2)}^2(E)@>{\theta\wedge}>>\cdots
\end{CD}}
\end{equation}
\begin{equation} {\begin{CD}
\{\mathcal{A}_{(2)}^\centerdot(E),
D''_u\}=\mathcal{A}_{(2)}^0(E)@>{D''_u}>>
\mathcal{A}_{(2)}^1(E)@>{D''_u}>>\mathcal{A}_{(2)}^2(E)@>{D''_u}>>\cdots
\end{CD}}
\end{equation}
where the sheaves $\mathcal{A}_{(2)}^i(E)$ are defined as follows:
\\
{\it for an open subset $U$ of $\overline M$,
$\mathcal{A}^i_{(2)}(E)(U)$ is defined as the set of $E$-valued
$i$-forms $\eta$ on $U\cap M$ with measurable coefficients and
measurable exterior derivative $D''_u\eta$, such that $\eta$ and
$D''_u\eta$ have finite $L^2$-norm};
\\
and the sheaves $\Omega_{(2)}^i(E)$ are defined as a sub-sheaves
of $\mathcal{A}_{(2)}^i(E)$ germs of which are local holomorphic
forms.
\\
\\
Define the hypercohomology of the complex
$\{\Omega_{(2)}^\centerdot, \theta\wedge\}$ as the {\it
$L^2$-Dolbeault cohomology}, denoted by $H^*_{{\text{Dol}},(2)}(M,
E)$. Since the complex $\{\mathcal{A}_{(2)}^\centerdot(E),
D''_u\}$ is a complex of fine sheaves, its hypercohomology is
computed by the cohomology of the corresponding complex of global
sections $\{A_{(2)}^\centerdot(E), D''_u\}$; call this cohomology
the {\it $L^2$-Higgs cohomology}, denoted by
$H^*_{{\text{Higgs}},(2)}(M, E)$. Then, one has the following
\\
\\
{\bf Conjecture 2.} {\it The $L^2$-Dolbeault cohomology
$H^*_{{\text{Dol}},(2)}(M, E)$ is isomorphic to the $L^2$-Higgs
cohomology $H_{{\text{Higgs}},(2)}^*(M, E)$; equivalently, the
complex $\{\Omega_{(2)}^\centerdot, \theta\wedge\}$ is
quasi-isomorphic to the complex
$\{\mathcal{A}_{(2)}^\centerdot(E), D''_u\}$.}
\\
\\
As in the compact case, using the theory of harmonic forms, one
can easily show that the $L^2$-de Rham cohomology
$H_{{\text{DR}},(2)}^*(M, \mathbb{V})$ is isomorphic to the
$L^2$-Higgs cohomology $H_{{\text{Higgs}},(2)}^*(M, E)$. To this
end, let us first recall the local K\"ahler identity for harmonic
bundles.
\\
\\
Using the previous notation, one has the first order K\"ahler
identities \cite{si1}
\begin{eqnarray*}
&&(D'_u)^*=\sqrt{-1}[\Lambda, D''_u],
~~~(D''_u)^*=-\sqrt{-1}[\Lambda,
D'_u]\\
&&(D^c_u)^*=-\sqrt{-1}[\Lambda, D], ~~~(D)^*=\sqrt{-1}[\Lambda,
D^c_u],
\end{eqnarray*}
where $^*$ represents the adjoint of the respective operator by
using the metrics $\omega$ and $u$. On the other hand, set
\[
\Delta=DD^*+D^*D, ~~~\Delta''=D''_u(D''_u)^*+(D_u'')^*D''_u.
\]
Using the above first order identities, one then has
\[
\Delta=2\Delta''.
\]
This shows that spaces of $\Delta$-harmonic forms valued in the
local system $\mathbb{V}$ can be identified with that of
$\Delta''$-harmonic forms valued in the Higgs bundle $E$.
\\
On the other hand, by $L^2$-theory, the $L^2$-de Rham cohomology
$H_{{\text{DR}},(2)}^*(M, \mathbb{V})$ can be represented by $L^2$
$\Delta$-harmonic forms valued in $\mathbb{V}$ and the $L^2$-Higgs
cohomology $H_{{\text{Higgs}},(2)}^*(M, E)$ can be represented by
$L^2$ $\Delta''$-harmonic forms valued in $E$. So,
$H_{{\text{DR}},(2)}^*(M, \mathbb{V})\cong
H_{{\text{Higgs}},(2)}^*(M, E)$.

\begin{rem}
The intersection cohomology is a topological notion \cite{gm1,
gm2}; but from the definition of the $L^2$-de Rham cohomology (and
the $L^2$-Higgs cohomology), we know that it depends on the
metrics $\omega$ and $u$, so that it is an analytical notion. From
the next two sections, we will see that the $L^2$-Dolbeault
cohomology is a much more algebraic notion, which is closely
related to the properties at the divisor of the representation
$\rho$. Thus, Conj. 1 and 2 establish the relations between these
different notions. When the harmonic bundle comes from a variation
of Hodge structure, these conjecture was first suggested by
Deligne.
\end{rem}

\noindent The conjectures were first proved by S. Zucker \cite{z}
when $M$ is a Riemann surface and the harmonic bundle comes from a
variation of Hodge structure. It should be pointed out that in the
Riemann surface's case, the intersection cohomology
$H^*_{\text{int}}(\overline M, j_*\mathbb{V})$ is isomorphic to
the sheaf cohomology $H^*(\overline M, j_*\mathbb{V})$, but this
is not in general valid in the higher dimensional case.
\\
\\
In the case that $M$ is a higher dimensional quasi-compact
K\"ahler manifolds and $\mathbb{V}$ is a variation of Hodge
structure, Conjecture 1 was proved independently by
Cattani-Kaplan-Schmid \cite{cks1} and Kashiwara-Kawai \cite{kk}.
We also should mention the works of  Looijenga \cite{lo} and
Saper-Stern \cite{ss}; they showed independently that the
$L^2$-cohomology is isomorphic to the intersection cohomology when
the base manifold is an arithmetic variety but the representation
is induced from a linear representation of the corresponding Lie
group.
\\
\\
In the note, we survey our recent study about Conjecture 1 and 2
\cite{jyz, jyz1}. In Section 2, we outline the work when the
coefficient is a VHS; Section 3 outlines the work when the base
manifold is a noncompact algebraic curve and the representation is
unipotent. A natural problem is what the situation is when we
don't assume that the representation $\rho$ is unipotent.

\section{Cohomologies valued in VHSs}
In this section, we consider cohomologies with coefficients in a
variation of Hodge structure and outline the results in
\cite{jyz}.
\\
\\
Let $H_{\mathbb{C}} =
H_{\mathbb{R}}\otimes_{\mathbb{R}}{\mathbb{C}}$, $k$ a positive
integer, and $S$ a real nondegenarate bilinear form on
$H_{\mathbb{C}}$ satisfying $S(u, v)=(-1)^kS(v, u)$. A {\it
polarized Hodge structure} of weight $k$ is a decomposition
$H_{\mathbb{C}}=\sum_{p+q=k}H^{p,q}$ satisfying
$H^{p,q}=\overline{H^{q,p}}$ and
\begin{eqnarray*}
S(H^{p,q}, H^{r,s}) &=& 0, ~~ {\text {unless}} ~ p=s, q=r,  \\
S({\mathcal C}v, {\overline v}) &>& 0, ~~ {\text {if}} ~ v\in
H^{p,q}, v \neq 0,
\end{eqnarray*}
where ${\mathcal C}$ is the Weil operator defined by ${\mathcal
C}v = {\sqrt{-1}}^{p-q}v, v\in H^{p,q}$. Call
$\dim_{\mathbb{C}}H^{p,q}$ the {\it Hodge number}, denoted by
$h^{p,q}$. Let $\bf D$ be the classifying space of polarized Hodge
structures of weight $k$, which is a homogeneous complex manifold.
Denote by $G_{\mathbb{R}}$ the group of holomorphic automorphisms
of $\bf D$ (more precisely, the group preserving $S$ and
$H_{\mathbb{R}}$), by $G_{\mathbb{C}}$ the complexification of
$G_{\mathbb{R}}$; denote by $\mathfrak{g}_0$ (resp.
$\mathfrak{g}$) the Lie algebra of $G_{\mathbb{R}}$ (resp.
$G_{\mathbb{C}}$).
\\
\\
One can consider a family of polarized Hodge structures
parameterized by a complex manifold $M$; this is the notion of
{\it (rational) variation of polarized Hodge structures} (in
brief, VHS). For the precise definition, one can refer to e.g.
\cite{sc}, p. 220. As in the introduction, we always assume that
$M={\overline M}-D$ is a quasi-compact K\"ahler manifold with a
Poincar\'e-like metric $\omega$, $D$ being a normal crossing
divisor, $j: M\to {\overline M}$ the inclusion map.
\\
\\
Let $\{M, {\bf H}_{\mathbb{Z}}\subset{\bf H}_{\mathbb{C}}, \{{\bf
F}^p\}_{p=0}^k, \nabla =\nabla^{1,0}+\nabla^{0,1}, {\bf S}\}$ be a
VHS, where ${\bf H}_{\mathbb{C}}$ is a flat bundle under a flat
connection $\nabla$, ${\bf H}_{\mathbb{Z}}$ is a flat lattice,
$\{{\bf F}^p\}_{p=0}^k$ is a Hodge filtration, and $\bf S$ is a
polarization. The polarization $\bf S$ defines a Hermitian metric
\[
h(\cdot, \cdot)={\bf S}(\mathcal{C}\cdot, \overline{\cdot}),
\]
called {\it Hodge metric}, where the bar is the conjugation with
respect to ${\bf H}_{\mathbb{Z}}$. By $\Vert\cdot\Vert$ denote the
corresponding norm. Take the successive quotients ${\bf F}^p/{\bf
F}^{p-1}$, called {\it Hodge bundles}, denoted by ${\bf E}^p$.
Accordingly, we have the induced $\cal{O}$-linear map $\theta^p:
{\bf E}^p\to{\bf E}^{p-1}$ of $\nabla^{1,0}$. Set ${\bf
E}=\oplus{\bf E}^p$ and $\theta=\oplus\theta^p$. Since
$(\nabla^{1,0})^2=0$, $\theta\wedge\theta=0$. So, the pair $({\bf
E}, \theta)$ is a Higgs bundle. One can also assign to the
variation a ($\rho$-equivariant) period mapping $\tilde\phi
:\tilde{M}\to {\bf D}$ from the universal covering $\tilde M$ of
$M$ to $\bf D$, which is holomorphic; here, $\rho$ is the induced
linear representation of $\pi_1(M)$ by $\nabla$.
\\
\\
Using the Poincar\'e-like metric $\omega$ on $M$ and the Hodge
metric $h$ on ${\bf H}_{\mathbb{C}}$, we define ${\Omega}^r({\bf
H}_{\mathbb{C}})_{(2)}$ to be the sheaf of germs of local $L^2$
holomorphic $r$-forms valued in $j_*{\bf H}_{\mathbb{C}}$ on
$\overline M$, where $j$ is the inclusion map of $M$ into
$\overline M$. According to the Hodge filtration $\{{\bf F}^p\}$
of the variation, one can then construct a filtration of
$\Omega^r({\bf H}_{\mathbb{C}})_{(2)}$, denoted by
$F^{p}\Omega^r({\bf H}_{\mathbb{C}})_{(2)}$, which is the sheaf of
germs of local $L^2$-holomorphic $r$-forms on $\overline M$ with
values in $j_*{\bf F}^{p-r}$. On the other hand, the quotient
bundle ${\bf E}^{p-r}={\bf F}^{p-r}/{\bf F}^{p-r+1}$ inherits a
quotient norm from ${\bf F}^{p-r}$; using this quotient norm, we
define $(\Omega^r\otimes {\bf E}^{p-r})_{(2)}$ as the sheaves of
germs of local $L^2$-holomorphic forms on $\overline M$ with
values in $j_*{\bf E}^{p-r}$, which, as seen in the following
proposition 1, is actually $(j_*{\bf
E}^{p-r}\otimes\Omega^r_{\overline M}(\log D))_{(2)}$, the sheaf
of germs of local $L^2$-sections in $j_*{\bf
E}^{p-r}\otimes\Omega^r_{\overline M}(\log D)$. Then, one can show
the following
\begin{thm}
The projection ${\bf F}^{p-r}\to {\bf E}^{p-r} (p\ge r)$ induces a
natural projection
\[
F^{p}\Omega^r({\bf H}_{\mathbb{C}})_{(2)}\to (\Omega^r\otimes {\bf
E}^{p-r})_{(2)};
\]
and the sequence
\begin{eqnarray*}
0\to F^{p+1}\Omega^r({\bf H}_{\mathbb{C}})_{(2)}\hookrightarrow
F^{p}\Omega^r({\bf H}_{\mathbb{C}})_{(2)}\to (\Omega^r\otimes {\bf
E}^{p-r})_{(2)}\to 0
\end{eqnarray*}
is exact, where $\hookrightarrow$ is the inclusion map.
\end{thm}
\begin{rem}
Theorem 1 was first proved by Zucker in the case of Riemann
surfaces \cite{z}.
\end{rem}
Since the following two results are local, w.l.o.g., we can assume
$M=(\triangle^*)^m$, $\triangle^*$ being the punctured disk.
$\pi_1(M)$ is then a free abelian group generated by $m$ elements
$\sigma_1, \sigma_2, \cdots , \sigma_m$, $\sigma_i$ corresponding
to the counter-clockwise path around $0$ of the $i$-th component
of ${(\triangle}^*)^m$. Denote by $\gamma_i$ the image of
$\sigma_i$ under $\rho$, which is possibly trivial and (if
nontrivial) referred to as the $i$-th {\it monodromy
transformation} of the variation. By a result of Borel \cite{sc},
up to some finite lifting, {\it throughout this paper we assume
that each $\gamma_i$ is unipotent}. Set $N_i=\log\gamma_i$ with
the elements on the diagonal being zero.
\begin{cor}
$\theta$ and hence $\nabla^{1,0}$ have the following asymptotic
behavior at $(0, 0, \cdots, 0)$
\[
\sum_{i=1}^nN_i{\frac{dz_i}{z_i}};
\]
and are $L^2$-bounded.
\end{cor}
For later applications, one needs some explicit expressions of the
sheaves $\Omega^r({\bf H_{\mathbb{C}}})_{(2)}$. This is done by
the following proposition 1, which says that the sheaf
$\Omega^r({\bf H_{\mathbb{C}}})_{(2)}$ can be defined
algebraically, just using the logarithmic monodromies $N_1,
\cdots, N_m$ and the corresponding weight filtrations, and lies in
$j_*{\bf H}_{\mathbb{C}}\otimes\Omega^r_{\overline
M}({\text{log}}D)$, $r=1, \cdots, m$. As a consequence of this
fact together with the asymptotic behavior of $\theta$ and the
theorem 1, we obtain the $L^2$ holomorphic Dolbeault complex on
$\overline M$
\[
(*)~~~~~~~~~{\bf E}_{(2)}\stackrel{\theta}{\to}({\bf
E}\otimes\Omega^1_{\overline M}(\log D))_{(2)}
\stackrel{\theta}{\to}({\bf E}\otimes\Omega^2_{\overline M}(\log
D))_{(2)}\stackrel{\theta}{\to}\cdots,
\]
which is furthermore independent of the Poincar\'e-like metric and
the Hodge metric. The key point of the proof is the norm estimates
\cite{cks} of the Hodge metric near the singularity.

\begin{prop} The sheaves $\Omega^0({\bf H_{\mathbb{C}}})_{(2)},
\Omega^1({\bf H_{\mathbb{C}}})_{(2)}, \cdots, \Omega^n({\bf
H_{\mathbb{C}}})_{(2)}$ are determined by the monodromies $N_1,
\cdots, N_n$ and can be expressed in terms of the weight
filtrations $\{{\bf W}^{\bf j}_*\}, j=1, \cdots, n$ on the domain
$D_{\epsilon} = \{(t_1, \cdots,
t_n)~|~{\frac{\log|t_1|}{\log|t_2|}}>\epsilon, \cdots,
{\frac{\log|t_{n-1}|}{\log|t_n|}}>\epsilon,
-\log|t_n|>\epsilon\}$, $\epsilon>0$. More precisely, if
considering the case of $\dim M =2$, one has the following
explicit formulae, {\small{\begin{eqnarray*} &&\Omega^{0}({\bf
H}_{\mathbb{C}})_{(2)} = t_1t_2\overline{\bf H}_{\mathbb{C}} +
t_1\bigcup_{l_2-l_1\le 0} \overline{\bf W}_{l_1l_2} +
t_2\overline{\bf W}^{\bf 1}_{0} +
\bigcup_{l_1\le 0}^{l_2\le l_1}\overline{\bf W}_{l_1l_2}, \\
&&\Omega^{1}({\bf H}_{\mathbb{C}})_{(2)} =
{\frac{dt_1}{t_1}}\otimes \big(t_1t_2\overline{\bf H}_{\mathbb{C}}
+ t_1\bigcup_{l_2-l_1\le 0}\overline{\bf W}_{l_1l_2} +
t_2\overline{\bf W}^{\bf 1}_{-2} + \bigcup_{l_1\le -2}^{l_2\le
l_1}\overline{\bf W}_{l_1l_2}\big) +     \\
&&+ {\frac{dt_2}{t_2}}\otimes \big(t_1t_2\overline{\bf
H}_{\mathbb{C}} + t_1\bigcup_{l_2\le l_1-2}\overline{\bf
W}_{l_1l_2} + t_2\overline{\bf W}^{\bf 1}_{0} + \bigcup_{l_1\le
0}^{l_2\le l_1-2}\overline{\bf W}_{l_1l_2}\big),  \\
&&\Omega^{2}({\bf H}_{\mathbb{C}})_{(2)} =
{\frac{dt_1}{t_1}}\wedge{\frac{dt_2}{t_2}}\otimes\big(
t_1t_2\overline{\bf H}_{\mathbb{C}} + t_1\bigcup_{l_2\le
l_1-2}\overline{\bf W}_{l_1l_2} + t_2\overline{\bf W}^{\bf 1}_{-2}
+ \bigcup_{l_1\le -2}^{l_2\le l_1-2}\overline{\bf W}_{l_1l_2}
\big),
\end{eqnarray*}}}
where $\overline{\bf W}_{l_1l_2}=\overline{\bf W}^{\bf
1}_{l_1}\cap \overline{\bf W}^{\bf 2}_{l_2}$ and ${\tilde
v}\in\overline{\bf W}^{\bf 1}_{l_1}\cap \overline{\bf W}^{\bf
2}_{l_2}$, ${l_2-l_1\le 0}$, implies that $v$ has nontrivial
projections to both $Gr^{{\bf W}^{\bf 1}_*}_{l_1}$ and $Gr^{{\bf
W}^{\bf 2}_*}_{l_2}$.
\end{prop}
In the above, we carefully analyze the $L^2$ holomorphic Dolbeault
complex $\{({\bf E}\otimes\Omega^r_{\overline M}(\log D))_{(2)},
\theta\}$ of a VHS on $\overline M$. For uniformity of notations,
we from now on denote $({\bf E}^{p-r}\otimes\Omega^r_{\overline
M}(\log D))_{(2)}$ by ${\text Gr}_{F}^p{\Omega}^{r}({\bf
H}_{\mathbb{C}})_{(2)}$; it is easy to see that it is actually a
piece of $({\bf E}\otimes\Omega^r_{\overline M}(\log D))_{(2)}$.
\\
\\
By the infinitesimal period relation of ${\nabla}^{1,0}$ and the
$L^2$-boundedness of $\theta$, we have
\[
\theta({\text Gr}_{F}^p{\Omega}^{r}({\bf
H}_{\mathbb{C}})_{(2)})\subset {\text Gr}_{F}^p{\Omega}^{r+1}({\bf
H}_{\mathbb{C}})_{(2)}.
\]
(More precisely, we should restrict to a piece of $\theta$.) Thus,
we obtain a holomorphic Dolbeault subcomplex of $\{({\bf
E}\otimes\Omega^{\cdot}_{\overline M}(\log D))_{(2)}, \theta\}$ on
$\overline M$
\[
\{{\text Gr}_{F}^p{\Omega}^{\cdot}({\bf H}_{\mathbb{C}})_{(2)},
\theta\}.
\]
Let $F^pA^k({\bf H}_{\mathbb{C}})_{(2)} =
\oplus_{r+s=k}(A^{r,s}\otimes {\bf F}^{p-r})_{(2)}$ and $D''_1 =
{\overline{\partial}} + {\nabla}^{1,0}$ for $p\ge 0$. Here,
$A^{r,s}$ is the sheaf of germs of local forms of type $(r, s)$
(not necessarily smooth) on $\overline M$ and $(A^{r,s}\otimes
{\bf F}^{p-r})_{(2)}$ is the sheaf of germs of local $L^2$ ${\bf
F}^{p-r}$-valued forms $\phi$ of type $(r,s)$ on $\overline M$ for
which ${\overline\partial}\phi$ are $L^2$ in the weak sense, where
the action of $\overline{\partial}$ is defined as follows: Let
$\phi$ be a form of type $(r, s)$ and $v$ a holomorphic section of
${\bf F}^p$, then $ \overline\partial(\phi\otimes
v)=\overline\partial\phi\otimes v$.
\\
\\
By the boundedness of $\nabla^{1,0}$, $D''_1(F^pA^k({\bf
H}_{\mathbb{C}})_{(2)})\subset F^pA^{k+1}({\bf
H}_{\mathbb{C}})_{(2)}$. Using the Hodge filtration $\{{\bf
F}^p\}$, take the successive quotients
\[
[{\text Gr}_F^pA^k({\bf H}_{\mathbb{C}})]_{(2)} := F^pA^k({\bf
H}_{\mathbb{C}})_{(2)}/F^{p+1}A^k({\bf H}_{\mathbb{C}})_{(2)},
\]
which, by the theorem 1 (which is also true for the differentiable
case by the proof), can be identified with
$\oplus_{r+s=k}(A^{r,s}\otimes {\bf E}^{p-r})_{(2)}$. Here,
$(A^{r,s}\otimes {\bf E}^{p-r})_{(2)}$ is the sheaf on $\overline
M$ of germs of local $L^2$ ${\bf E}^{p-r}$-valued forms $\phi$ of
type $(r,s)$ for which ${\overline \partial}\phi$ are $L^2$ in the
weak sense. Denote the induced map of $D_1''$ by $D''$, which is
actually ${\overline{\partial}} + \theta$ and satisfies $(D'')=0$.
We now obtain a complex of fine sheaves on $\overline M$
\[
\{[{\text Gr}_F^pA^{\cdot}({\bf H}_{\mathbb{C}})]_{(2)}, D''\},
\]
for $p\ge 0$ and the holomorphic Dolbeault subcomplex $\{{\text
Gr}_{F}^p{\Omega}^{\cdot}({\bf H}_{\mathbb{C}})_{(2)}, \theta\}$
is its subcomplex. Similar to the complex $\{{\text
Gr}_{F}^p{\Omega}^{\cdot}({\bf H}_{\mathbb{C}})_{(2)}, \theta\}$,
$([{\text Gr}_F^pA^{\cdot}({\bf H}_{\mathbb{C}})]_{(2)}, D'')$ can
also be considered as a piece of a larger complex of fine sheaves:
Denote by $[{\text Gr}_F^{*}A^{k}({\bf H}_{\mathbb{C}})]_{(2)}$
the sheaf of germs of local $L^2$ $k$-forms $\phi$ (not
necessarily smooth) with values in $\bf E$ on $\overline M$, for
which $\overline\partial\phi$ are also $L^2$ in the weak sense. It
is clear that $\{[{\text Gr}_F^{*}A^{\cdot}({\bf
H}_{\mathbb{C}})]_{(2)}, D''\}$ is a complex of fine sheaves and
$\{[{\text Gr}_F^pA^{\cdot}({\bf H}_{\mathbb{C}})]_{(2)}, D''\}$
is a piece of it. Then we can show the following
\begin{thm}
The holomorphic Dolbeault complex $\{{\text
Gr}_{F}^p{\Omega}^{\cdot}({\bf H}_{\mathbb{C}})_{(2)}, \theta\}$
is quasi-isomorphic to the complex $\{[{\text
Gr}_F^pA^{\cdot}({\bf H}_{\mathbb{C}})]_{(2)}, D''\}$ under the
inclusion map for $p\ge 0$.
\end{thm}
Theorem 2 then implies Conjecture 2 in the case of variation of
Hodge structure.

\section{Cohomologies valued in harmonic bundles on non-compact
curves} In this section, we consider cohomologies valued in
harmonic bundles on non-compact curves. In order to define the
$L^2$-cohomology, we need a suitable harmonic metric, which
reflects the geometry of the representation $\rho$. Throughout
this section, we assume that $\rho$ is unipotent near the divisor.

\subsection{Harmonic metrics in the case of noncompact curves}
From now on, we assume that $M$ is a noncompact curve, i.e. a
compact Riemannian surface deleting finitely many points, and
change the symbol $M$ into $S$, the compactification of which is
denoted by $\overline S$; $S={\overline S}\setminus\{p_1, \cdots,
p_s\}$ and $j: S\to\overline S$ is the inclusion map. Under such a
case together with the assumption that the representation $\rho$
is unipotent, one can easily get a $\rho$-equivariant harmonic map
of finite enegy and its behavior near the punctures $\{p_1,
\cdots, p_s\}$. It is worth pointing out that in \cite{jz}, the
authors considered these problems without the restriction of
dimension. Since there the target manifolds are very general, the
construction of the initial metrics are getting very complicated.
In the present case, since the representation is linear, we will
give an explicit construction though the basic idea is the same as
\cite{jz}; more importantly, we can see explicitly why the initial
map is of finite energy, and hence the analysis of \cite{jz} is
applicable. The idea of the construction for the initial map will
be used in the most general case, where representations need not
be unipotent and also the energy of maps need not be finite. Here,
we should also point out that the construction in \cite{jz1} is
also a very special one (cf. \cite{jyz2}); in \cite{y, y1}, we
consider some constructions with some artificial singularities. As
before, we take the Poincar\'e-like metric on the base manifold
$S$, and by $t$ denote local Euclidean complex coordinate.
\\
\\
Let $\rho: \pi_1(S)\to {GL}(n, \mathbb{C})$ be a semisimple linear
representation, and restrict $\rho$ to a neighborhood of $p_i$,
say a small punctured disk $\Delta^*$ around $p_i$, which we call
the {\it boundary representation} of $\rho$, denoted by $\rho_i$.
Call the representation $\rho$ {\it unipotent} if each
$\rho_i(\gamma)$ is a unipotent matrix, where $\gamma$ is a circle
around $p_i$. A result of Borel tells us that this is the case for
VHS up to a suitable lifting (cf. e.g. \cite{sc}, Lemma 4.5). Take
the logarithm of $\rho_i(\gamma)$, denoted by $N$, which is
upper-triangle and all the diagonal entries of which are $0$ under
a suitable basis.
\\
\\
We now proceed to construct an initial metric on $L_{\rho}$,
equivalently, an initial $\rho$-equivariant map from $\tilde S$ to
$Gl(n, \mathbb{C})/U(n)$. To this end, let us first give some
preliminary. Let $\mathcal{P}_n$ be the set of all positive
definite hermitian symmetric matrices of order $n$. $Gl(n,
\mathbb{C})$ acts transitively on $\mathcal{P}_n$ by
\[
g\circ H=:{g}H^t\bar{g}, ~H\in\mathcal{P}_n, g\in GL(n,
\mathbb{C}).
\]
Obviously, the action has the isotropic subgroup $U(n)$ at the
identity $I_n$. Thus $\mathcal{P}_n$ can be identified with the
coset space $Gl(n, \mathbb{C})/U(n)$, and can be uniquely endowed
an invariant metric{\footnote{In terms of matrices, such an
invariant metric can be defined as follows. At the identity $I_n$,
the tangent elements just are hermitian matrices; let $A, B$ be
such matrices, then the Riemannian inner product $<A,
B>_{\mathcal{P}_n}$ is defined by $tr(AB)$. In general, let $H\in
\mathcal{P}_n$, $A, B$ two tangent elements at $H$, then the
Riemannian inner product $<A, B>_{\mathcal{P}_n}$ is defined by
$tr(H^{-1}AH^{-1}B)$.} up to some constants. In particular, under
such a metric, the geodesics through the identity $I_n$ are of the
form $\exp(tA)$, $t\in \mathbb{R}$, $A$ being a hermitian matrix.
\\
\\
Let the Jordan normal form of $N$ have $p$ Jordan blocks, denoted
by $N_j, 1\le j\le p$. By the Jacobson-Morosov theorem, one can
expand each block $N_j$ into an $\mathfrak{sl}_2$-triple $\{Y_j,
N_j, N_j^-\}$, i.e., $[Y_j, N_j]=2N_j, [Y_j, N_j^-]=-2N_j^-$ and
$[N_j, N_j^-]=Y_j$; a theorem of Kostant tells us that such a
triple, up to conjugations, is unique (cf. e.g. \cite{ko}). Take
the Euclidean coordinate $t$ on $\Delta^*$ with $t(p_i)=0$ and
$t=re^{\sqrt{-1}x}$; also take the universal covering
$$
H_{\alpha}=\{z=x+\sqrt{-1}y~|~x\in \mathbb{R}, ~y>-\log\alpha\}
$$
of $\Delta^*$ with $y=-\log r$, for a positive number $\alpha<1$.
Corresponding to a (once and for all) fixed flat sections basis of
$L_\rho$, we construct the required Hermitian metric of $L_{\rho}$
on $\Delta^*$ as
\begin{equation}\label{metric}
h_i=\left(
\begin{array}{cccc}
M_1&~&0\\
~&\ddots&~\\
0&~&M_p
\end{array}
\right),
\end{equation}
where $ M_j=\exp (xN_j)\circ\exp(({\frac 1 2}\log|\log r|)Y_j).$
In the following, $h_i$ is considered as both a metric and the
above matrix under the fixed basis. Clearly, $h_i$, as a map from
$H_{\alpha}$ to $\mathcal{P}_n$, is $\rho_i$-equivariant when
changing $\log r$ into $y$. Geometrically, this is a geodesic
$\rho_i$-equivariant embedding of $\Delta^*$ into $\mathcal{P}_n$
which maps the puncture to the infinity of $\mathcal{P}_n$, as can
be explicitly seen when one embeds geodesically the upper-half
plane into the real hyperbolic $3$-space ${\bf H}^3$ by using the
matrix model of ${\bf H}^3$ (for more details, cf. \cite{jyz2}).
We also remark that each $h_i$ is actually harmonic. Now one can
easily extend the metrics $\{h_i\}$ to a global metric on
$L_\rho$, denoted by $h_0$. Naturally, the corresponding map is
$\rho$-equivariant.
\\
\\
We now want to show that each $h_i$, and hence $h_0$, is of finite
energy. For simplicity, we may assume here that the Jordan normal
form of $N$ has only one Jordan block. Let $\{Y, N, N^-\}$ be the
corresponding $\mathfrak{sl}_2$-triple. The semisimple element $Y$
can actually be described as follows. Canonically, $\mathbb{C}^n$
has a filtration
\begin{equation}\label{filtr}
0\subset W_{-(n-1)}\subset W_{-(n-3)}\subset\cdots\subset
W_{n-3}\subset W_{n-1}=\mathbb{C}^n,
\end{equation}
satisfying that $N(W_i)\subset W_{i-2}$, $Y$ preserves each $W_i$,
and all the quotients $W_i/W_{i-2}$ are $1$-dimensional. Then the
(induced) action of $Y$ on $W_i/W_{i-2}$ is multiplying by $i$.
Actually, one can also choose a basis $\{e_{-(n-1)}, e_{-(n-3)},
\cdots, e_{n-3}, e_{n-1}\}$ of $\mathbb{C}^n$, which is compatible
with the above filtration (i.e., $Ne_j=e_{j-2}$ and $\{e_j\}_{j\le
i}$ generates $W_i$) and satisfies $Ye_i=ie_i$. $\exp(({\frac 1
2}\log|\log r|)Y)$, under the above basis, can then be written as
\begin{equation}\label{metric1}
\left(
\begin{array}{ccccc}
|\log r|^{\frac{-(n-1)}{2}}&0&\cdots&0&0\\
0&|\log r|^{\frac{-(n-3)}{2}}&\cdots&0&0\\
\vdots&\vdots&\ddots&\vdots&\vdots\\
0&0&\cdots&|\log r|^{\frac{n-3
}{2}}&0\\
0&0&\cdots&0&|\log r|^{\frac{n-1}{2}}
\end{array}
\right).
\end{equation}
Then, using the invariant metric of $\mathcal{P}_n$, a simple
computation shows that the energy of $h_i$ satisfies
\[
E(h_i)=\int_{\Delta^*}|h_i^{-1}dh_i|^2*1\le C\int_{0}^\alpha |\log
r|^{-2}r^{-1}dr<\infty.
\]
\\
\\
As mentioned before, since the initial map $h_0$ is of finite
energy, the analysis of \cite{jz} works; in particular, one has
the following
\begin{prop}
Let $\rho: \pi_1(S)\to Gl(n, \mathbb{C})$ be a semisimple
representation all the boundary representations of which are
unipotent, $h_0$ be the $\rho-$equivariant map (initial metric)
constructed above. Then there exists a $\rho$-equivariant harmonic
map (harmonic metric) of finite energy
\[
h: {\tilde S}\to Gl(n, \mathbb{C})/U(n),
\]
which has the same asymptotic behavior as $h_0$ near the
punctures, where $\tilde S$ is the universal covering of $S$;
moreover the norm of the derivative $dh$ of $h$ satisfies, when
going down to $S$ and measured near the divisor with respect to
the Poincar\'e-like metric and the standard Riemannian symmetric
metric on $Gl(n, \mathbb{C})/U(n)$,
$$
|dh|^2\le C{|\log r|^2}
$$
for some constant $C>0$, where $r$ is the radial Euclidean
coordinate of $\Delta^*$.
\end{prop}
\begin{rem}
From (11) and the related argument there, one has some explicit
asymptotic estimates of the harmonic metric $h$, which are the
same as the case of VHS (cf. \cite{sc}).
\end{rem}

\subsection{The $L^2$-cohomology: The $L^2$-Poincar\'e Lemma}
For the direct image sheaf $j_*L_{\rho}$ of the local system
$L_{\rho}$ on $\overline S$, one has the \v{C}ech cohomology
$H^*({\overline S}, j_*L_{\rho})$; on the other hand, using the
Poincar\'e-like metric $\omega$ on $S$ and the harmonic metric $h$
on $L_{\rho}$, one can define a complex
$\{\mathcal{A}^{\cdot}_{(2)}(L_{\rho}), D\}$ of sheaves over
$\overline S$ as in the introduction and then $L^2$-de Rham
cohomology $H^*(\{\Gamma(\mathcal{A}^{\cdot}_{(2)}(L_{\rho})),
D\})$, denoted by $H^*_{{DR}, (2)}(\overline S, L_\rho)$.

\begin{thm}
There exists a natural identification
\[
H^*({\overline S}, j_*L_{\rho})\cong H^*_{{DR}, (2)}(\overline S,
L_\rho).
\]
\end{thm}
This is Conjecture 1 in the present situation.
\\
\\
Canonically, the proof of Theorem 1 is reduced to prove

\begin{thm} (The $L^2$-Poincar\'e lemma)
The complex $\{\mathcal{A}^{\cdot}_{(2)}(L_{\rho}), D\}$ is a
resolution of $j_*L_{\rho}$. This is equivalent to saying that

\noindent 1)
$j_*L_{\rho}=\{\eta\in\mathcal{A}^{0}_{(2)}(L_{\rho})~|~D\eta=0\}$;

\noindent 2) the differential $D$ satisfies the Poincar\'e lemma,
i.e., if an $i$-form $\eta$ in $\mathcal{A}^{i}_{(2)}(L_{\rho})$
is $D$-closed, then there exists an $i-1$-form $\sigma\in
\mathcal{A}^{i-1}_{(2)}(L_{\rho})$ satisfying $D\sigma = \eta$,
for $i=1, 2$.
\end{thm}
\begin{rem}
In the case of higher dimension, one needs to consider the
intersection cohomology $H^*_{\text{int}}(\overline S, j_*L_\rho)$
instead of the \v{C}ech cohomology $H^*({\overline S},
j_*L_{\rho})$.
\end{rem}

\subsection{The $L^2$-Higgs cohomology: The $L^2$-$\overline\partial$-Poincar\'e Lemma}
As seen in the introduction, the harmonic metric $h$ on $L_{\rho}$
induces the structure of a Higgs bundle on $L_{\rho}$: $(E,
D''=\overline\partial+\theta)$, satisfying $D=D'+D''$ with
$D'=\partial+\overline\theta$; moreover, the Hermitian connection
$\partial+\overline\partial$ w.r.t. the metric $h$ has bounded
curvature under the metric $h$ and the Poincar\'e-like metric so
that $E$ can be analytically extended to $\overline S$, denoted by
$j_*E$ as usual, which is especially coherent. Furthermore,
$\theta$ has a $\log$-singularity, i.e. $\theta\sim
{\frac{dt}{t}}N$. It is especially worth pointing out that by an
argument of Simpson (cf. \cite{si}), {\it the residue $N$ of
$\theta$ here coincides with the logarithmic monodromy $N$ in the
local system $L_{\rho}$}; so although under different bundle
structures, we have the same weight filtration under certain
suitable identification. Throughout this subsection, we consider
the Higgs bundle $(E, D''=\overline\partial+\theta)$ together with
the harmonic metric $h$, satisfying that the meromorphic sections
of $E$ at the punctures have the same estimates as the case of VHS
\cite{sc}, just forgetting that it comes from the local system
corresponding to the representation $\rho$.
\\
\\
As seen in the introduction, using the Poincar\'e-like metric
$\omega$ on $S$ and the harmonic matric $h$ on $(E, D'')$, one can
define a complex $\{\mathcal{A}^{\cdot}_{(2)}(E), D''\}$ of fine
sheaves on $\overline S$ and hence the corresponding cohomology
${H}^*(\{\Gamma(\mathcal{A}^{\cdot}_{(2)}(E)), D''\})$ of the
complex of global sections--the {\it $L^2$-Higgs cohomology} of
$\overline S$ valued in the Higgs bundle $(E,
D=\overline\partial+\theta)$, denoted by $H^*_{{\text{Higgs}},
(2)}({\overline S}, E)$.
\\
\\
The following lemma is a key.

\begin{lem}
$\theta$ is an $L^2$-bounded operator.
\end{lem}
{\it Proof}. The proof uses the norm estimates of $h$ in the sense
of Higgs bundles.
\\
\\
Based on this lemma, as seen in the introduction, one has a
sub-complex of $\{\mathcal{A}^{\cdot}_{(2)}(E), D''\}$---the
$L^2$-{\it holomorphic Dolbeault complex}
$\{\Omega^{\cdot}_{(2)}(E), \theta\}$ and the corresponding
hypercohomology $\mathbb{H}^*(\{\Omega^{\cdot}_{(2)}(E),
\theta\})$,  the $L^2$-{\it Dolbeault cohomology} of $\overline S$
valued in the Higgs bundle $(E, D=\overline\partial+\theta)$,
denoted by $H^*_{\text{Dol}, (2)}(\overline S, E)$. Then, we have
\begin{thm} ($L^2$-$\overline\partial$-Poincar\'e lemma)
The inclusion
$$
i: \{\Omega^{\cdot}_{(2)}(E), \theta\}\hookrightarrow
\{\mathcal{A}^{\cdot}_{(2)}(E), D''\}
$$
is a quasi-isomorphism; and hence one has
\[
H^*_{{\text{Higgs}}, (2)}({\overline S}, E)\cong H^*_{\text{Dol},
(2)}(\overline S, E).
\]
\end{thm}
The proof of the theorem is reduced to show the following lemma,
which was first proved by Zucker in the case of VHS (cf. \cite{z},
Proposition 6.4).
\begin{lem} Let $V$ be
a holomorphic line bundle on $\Delta^*$ with generating section
$\sigma$, and with a Hermitian metric satisfying
\[
\Vert\sigma\Vert^2\sim |\log r|^k, ~~ k\in\mathbb{Z}, ~k\neq 1.
\]
Then for every germ of an $L^2$ $(0, 1)$-form
$\phi=fd\overline{t}\otimes\sigma$ at the puncture, there exists
an $L^2$ section $u\otimes\sigma$ with $\overline\partial
u=fd\overline{t}$.
\end{lem}


\end{document}